\documentclass{article}
\usepackage{latexsym,amssymb,amsmath,amsfonts,epsf,rotating}
\usepackage[dvips]{epsfig}

\begin{document}

\title{Discrete Reanalysis of a New Model of \\ the Distribution of Twin Primes}
\author{P.F.~Kelly\footnote{patrick\_kelly@ndsu.nodak.edu} 
\ and Terry~Pilling\footnote{terry@mailaps.org} \\
Department of Physics \\
North Dakota State University \\
Fargo, ND, 58105-5566 \\
U.S.A.}

\date{June, 2001}
 
\maketitle

\begin{abstract}
Recently we have introduced a novel characterisation of the distribution
of twin primes that consists of three essential elements.
These are:  that the twins are most naturally viewed as a subsequence
of the primes themselves, that the likelihood of a particular prime
in sequence being the first element of a twin is akin to a fixed-probability
random event, and that this probability varies with $\pi_{1}$,
the count of primes up to this number, in a simple way.
Our initial studies made use of two unproven
assumptions:  that it was consistent to model this fundamentally discrete
system with a continuous probability density, and that the fact that an
upper-bound cut-off for prime separations exists could be consistently
ignored in the continuous analysis.
The success of the model served as {\it a posteriori} justification 
for these assumptions.
Here we perform the analysis using a discrete formalism -- not passing
to integrals -- and explicitly include a self-consistently defined cut-off.
In addition, we reformulate the model so as to minimise the input data needed.
\end{abstract}

\vspace{5mm}


\noindent
MCS:  11N05 (Primary) 11B05, 11A41 (Secondary)

\section{Introduction}

In two recent papers, an empirical model for the distribution of twin
primes was proposed~\cite{random} and some of its predictions were
developed~\cite{return}.
The foundation of this novel approach to the distribution of twins
is that the sequence of twins is most naturally studied in the context
of the primes rather than the natural numbers.

In~\cite{random}, empirical evidence was presented which strongly supports
the contention that within the set of prime numbers less than or equal to
some number $N$, twins (pairs of primes with arithmetic difference 2)
occur in the manner of fixed probability random events.
This fact lies at the heart of the model.
The probability is not constant however, rather it decreases with 
increasing $N$, in accord with the Hardy--Littlewood Conjecture.
The third essential ingredient of the model is that the manner in which
the probability changes can be expressed simply in terms of 
$\pi_{1}(N)$, the number of primes less than or equal to $N$.
The reader is referred to~\cite{random,return} for details.
Extensions of the empirical analysis by ourselves and
others~\cite{newwolf} verify the persistence of the model up to 
$N \sim 10^{13}$.

It was noted first by J.~Calvo~\cite{jorge}, and also independently by
M.~Wolf~\cite{newwolf}, that in the course of developing the 
model for the distribution of twins we have taken an essentially 
discrete system of {\it prime separations} ({\it i.e.,} the number
of singleton primes which occur between a pair of neighbouring twins),
and modelled it with a continuous distribution.
Another possible shortcoming of the model developed in~\cite{random}
is that when we normalise the distribution we integrate over all 
(continuous) prime separations from 0 (the most likely separation in
the distributional model, see~\cite{return}) to $\infty$.
It is formally impossible to take this limit at any finite $N$, or even
generally if the number of twins is infinite.
We did so in our analysis because we believed that the error introduced 
was quite small.
This was borne out by the apparent success of the model.

It is the aim of this paper to address the two concerns:  discrete 
analysis {\it versus} continuous, and taking into account the fact 
that for any $N$ there exists a maximum prime separation.
In the next section the model is briefly reviewed and reformulated.
In the following section we shall reanalyse with
sums rather than passing to the integral representation.
Second, we shall self-consistently set an upper bound to the separations
and incorporate its effects into our analysis.
As a test of consistency, the predictions for upper bounds will be
compared with the analysis of ``prime gaps'' in~\cite{return}.

\section{The Model}

The model that we consider is empirical in that it is derived from
a direct analysis of the distribution of twin primes less than
$2 \times 10^{11}$.
The essential feature which provides the key to the success of the
model is that the distribution of twins is considered in the context
of the primes alone rather than within the natural numbers.
The model is based upon the observation that twins less than
some number $N$ seem to occur
as fixed-probability random events in the sequence of primes.
That is, there is a characteristic distribution of {\it prime separations}
which may expressed in the form
\begin{equation}
{\cal P}(s, \pi_1) = A \, e^{-ms} \, .
\label{prob}
\end{equation}
Here, $\pi_1$ is the number of primes less than or equal to $N$,
$s$ is the {\it prime separation} (the number of unpaired {\it singleton}
primes occuring between two twins),
and $m$ is a decay parameter which is constant for a given $N$,
but varies with $\pi_{1}(N)$, while $A(\pi_{1}(N))$ is an overall 
constant which is fixed by normalisation.
${\cal P}(s, \pi_1)$ is the probability density that a given 
pair of twins in the sequence of primes up to $N$ has prime separation 
$s$.
When~(\ref{prob}) is assumed to be continuous and extending to infinity,
the condition that it be properly normalised, 
$\int_0^\infty {\cal P}(s, \pi_1) \, ds \equiv 1$,
constrains $A = m$.

We chose a representative sample of prime sequences
and determined the decay constants for each.
We began our analysis with $(5 \ 7)$, 
discarding the anomalous twin $(3 \ 5)$.
The variation of the decay parameters -- the slopes on a plot of 
log(frequency) {\it versus} separation -- is well-described by 
the following function:
\begin{equation}
- m ( \pi_1 ) = - \frac{m_0}{ \log(\pi_1) } \, ,
\label{em}
\end{equation}
where the constant, $m_0$, has been estimated to equal $1.321 \pm .008$
in \cite{random}.

\subsection{Reformulating the Model}

Our empirical model is founded upon a constructive procedure:
from an exact knowledge of the distribution of separations, we
determine the best-fit slope on a graph of $\log($frequency$)$ 
{\it vs.} separation, giving each datum {\it equal} weight.
The unfortunate aspect of this is that we are limited to the
data that we have collected.
In particular, various groups of researchers have {\it counted}
primes and twins but they do not appear to have kept detailed
counts of twin separations.
Further, we now believe that we understand better the difficulties
which led us to eschew characterisation of the behaviour of twins 
solely in terms of $\pi_{1}$ and $\pi_{2}$.

Let us now make a proper case for consideration of an empirical model
for the distribution of twins whose inputs are the counts $\pi_{1}(N)$
and $\pi_{2}(N)$.
Recalling \cite{random}, especially the formulae $(7)$ and $(15)$ 
and Figure 4, (reproduced here as Figure 1),
in which the statistical average separation $s_0$ is expressed as 
the number of singletons divided by the number of twins, and 
$\bar{s}$ is the reciprocal of the estimated slope, we write
\begin{equation}
s_0 = \frac{\pi_1 - 2 \pi_2}{\pi_2}
\quad \mbox{ and } \quad
\bar{s} = \frac{1}{m}  \, .
\label{s_average}
\end{equation}
It was readily apparent that the slopes $m$, determined as described
above, and quantities $1/s_0$ as in (\ref{s_average}), closely correspond
at large values of $N$.
It is equally apparent that they differ significantly at smaller values.
This is the region where the relatively strong enhancement of the few 
``large'' separation events had the greatest effect in reducing the 
magnitude of the computed slopes.
This in turn enabled the success of our simple and straightforward
empirical model for the variation of the slope with $\pi_1$.

{\begin{figure}[htb]
\begin{center}
\begin{turn}{-90}
\leavevmode
\epsfxsize=3.5in
\epsfysize=4.5in
\epsfbox{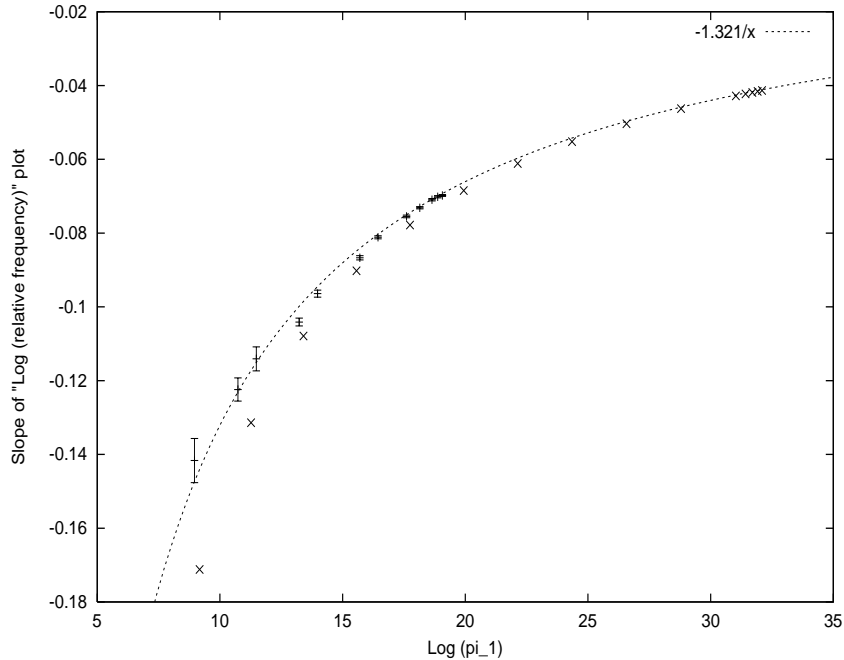}
\end{turn}
\end{center}
\caption{Synthetic slopes ($1\big/s_0$) using Nicely's data marked with
$\times$'s, computed slopes from the actual spectrum of prime separations
with error bars, and our empirical fit.  All {\it vs.}~$\log(\pi_1)$.
This is Figure 4 in \cite{random}.}
\label{graph1}
\end{figure}}

We set out to understand better the behaviour at the low end of the curve 
in Figure~\ref{graph1}, by reconsideration of $s_0$ as a function
of $\log(\pi_1)$.
In Figure~\ref{graph2} the $s_0$ derived from Nicely's 
data\footnote{Nicely's data~\cite{nicelydata} consist of values of $N$,
$\pi_{1}(N)$, and $\pi_{2}(N)$.  We have adjusted the $\pi_{1}$ 
to discard the singleton primes which appear after the last twin less 
than $N$.}
appear to follow very closely along a straight line with slope 
$0.7918 \pm 0.0007$ and $y$-intercept $-1.194 \pm 0.018$.
The negative value for the $y$-intercept, implying a positive value
for the $x$-intercept initially appeared to us to be pathological
and prevented us from arriving at a simple empirical characterisation 
for the variation of the Nicely data.
We now argue that this pathology is relatively benign,
as the $x$-intercept has such a small value, here $\sim 1.5$, that 
$\pi_{1}( N ) \simeq \exp(1.5)$ is less than 5, and thus 
the value of $N$ to which it corresponds is less than $20$.
We have no expectation that our statistical model can produce meaningful
results for short sequences of primes, and
so this value for the $x$-intercept is truly and completely extraneous%
\footnote{See the Conclusion for further comments about the 
accuracy of our linear fit and the values of the intercepts.}.

{\begin{figure}[htb]
\begin{center}
\begin{turn}{-90}
\leavevmode
\epsfxsize=3.5in
\epsfysize=4.5in
\epsfbox{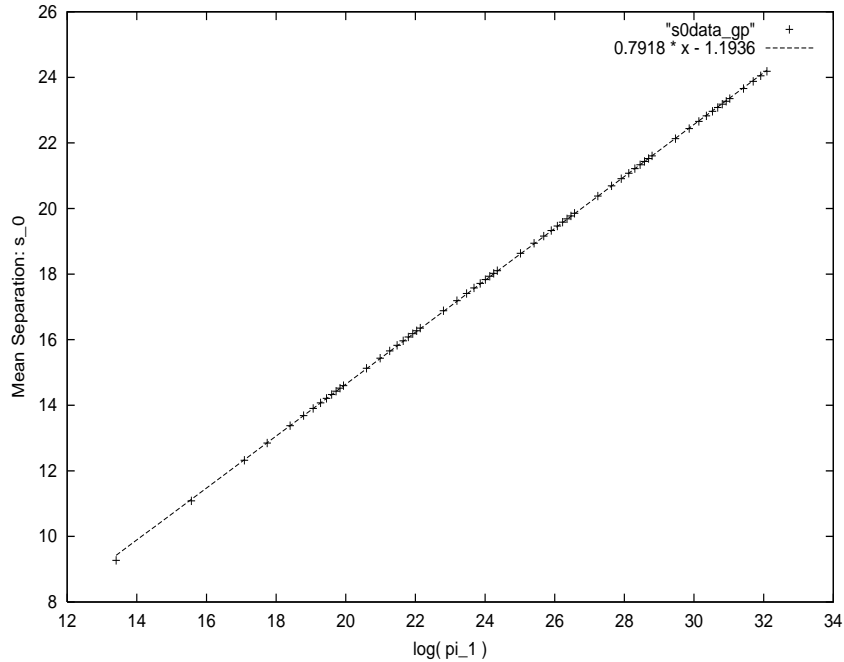}
\end{turn}
\end{center}
\caption{Statistical average prime separations, $s_0$, using 
Nicely's data and our linear fit {\it vs.}~$\log(\pi_1)$.}
\label{graph2}
\end{figure}}
 
Thus, for the purposes of this paper, we have reformulated the fundamentally
empirical model that we proposed in \cite{random,return} in such a way
that we characterise the probability distribution for the twin separations
not in terms of the ``decay constant'' $m$ and its proper variation with
$N$ as before, but instead in terms of the ``mean separation'' $\bar{s}$
and its concommitant proper variation.
It is most likely the case that this reformulation is best suited for
further development as it relies exclusively on data obtained solely
by counting primes and twins.
Henceforth we will rewrite (\ref{prob}) as 
\begin{equation}
{\cal P}(s, \pi_1) = A \, e^{-s/\bar{s}} \, .
\label{sprob}
\end{equation}
Our empirical model (Figure~\ref{graph2}) strongly suggests that 
$s_{0}$ varies with $\pi_{1}$ as
\begin{equation}
s_{0} ( \pi_1 ) = {\cal S}_{1} \log( \pi_{1} ) + {\cal S}_{0} \, ,
\label{snaughtempirical}
\end{equation}
and we've written ${\cal S}_{1}$ and ${\cal S}_{0}$ for the constants whose
values are empirically determined to be $0.7918 \pm 0.0007$, and
$-1.194 \pm 0.018$ respectively.
We further insist that $s_{0} > 0$ as discussed above.
The relation that exists between $\bar{s}$ and $s_{0}$ is revealed
in the next section.

\section{Discrete Analysis}

Normalising the probability distribution for the occurrence of 
prime separations yields
\begin{equation}
1 = \sum_{s = 0}^{L} {\cal P}(s, \pi_{1}) 
= \sum_{s = 0}^{L} A \, e^{- s / \bar{s}} \, ,
\label{norm}
\end{equation}
where we have interpreted the ${\cal P}( s , \pi_{1} )$ as relative 
frequencies rather than absolute counts (in which case the {\it lhs}
of (\ref{norm}) would equal $\pi_{2} - 2$, the total number of 
twins less two\footnote{ We disregard the twin $(3\ 5)$ and the 
prime separations are intervals between neighbouring twins.}
and have inserted as an {\it Ansatz} the empirical relation (\ref{sprob}).
$L$ denotes the maximum prime separation, thereby truncating the sum.
{\it A priori} $L$ is not specified and a value must be assumed 
or self-consistently derived.

A second relation among $A$, $\bar{s}$ and $L$ is formed by
consideration of the frequency-weighted average prime separation,
{\it viz.}
\begin{equation}
s_0 = \frac{\pi_{1} - 2\pi_{2}}{\pi_{2}}
= \sum_{s = 0}^{L} s \, {\cal P}(s,\pi_{1})
= \sum_{s = 0}^{L} s \, A \, e^{- s / \bar{s}} \, .
\label{snaught}
\end{equation}
In fact, a note of caution is required here.
Our empirical analysis discards the first twin and considers separations,
so ``$\pi_2$'' should be replaced by $\pi_{2} - 2$, and to be consistent
``$\pi_1$'' should be $\pi_{1} -2$ since we pass over the primes $2$ and $3$.
Incorporating these minor offsets into (\ref{snaught}) results in the
$(\pi_1 - 2 \pi_2 + 2)\big/(\pi_2 - 2)$ which is only slightly different
from the simpler and more straightforward expression that we use.

A third relation among the parameters comes from assigning a cut-off
for the probability distribution for prime separations.
The {\it Ansatz} (\ref{sprob}) has no such cut-off built into it,
although one might try to use the ``scale'' set by $\bar{s}$ to 
establish one by {\it fiat} (say, $L = 20 \times \bar{s}$).
Instead, we shall adopt the general method utilised in \cite{return},
and set a minimum probability threshold with the introduction of a 
so-called risk factor $f$.
In effect $f$ prime separation events with prime separations greater than
$L$ are be expected to occur in the context of the probabilistic 
model.
Then we may write
\begin{equation}
\frac{f}{\pi_{2}} = \sum_{s = L +1}^{\infty} {\cal P}(s,\pi_{1}) \, ,
\label{cut-off}
\end{equation}
providing a self-consistent cut-off value for the sum over separations.

Thus we have a set of three equations with three known quantities,
the counts of primes and twins, $\pi_{1}$ and $\pi_{2}$ respectively,
for a given $N$ and the risk factor $f$, and three parameters
to be determined: $A$, $\bar{s}$, and $L$.
We now proceed to solve these equations in two distinct instances:
the first in which we formally set $f = 0$ in which case $L = \infty$,
and the second, where $f$ is specified (non-zero).

\subsection{With $f = 0$ and $L = \infty$}

With $f = 0$ and $L = \infty$ equation (\ref{cut-off}) has no content
and (\ref{norm}) and (\ref{snaught}) can be solved exactly.
The solutions are
\begin{equation}
\frac{1}{\bar{s}} = \log \left[ 1 + \frac{1}{s_0} \right] \ ,
\label{soln1s}
\end{equation}
and
\begin{equation}
A = \frac{1}{ 1 + s_{0} } \ .
\label{soln1A}
\end{equation}
We note three things about this solution.
The first is that $\bar{s}$ and $A$ depend on $\pi_{1}$ and $\pi_{2}$
and hence $N$ {\it implicitly} via $s_0$.
The second is that with increasing $N$ the density of twins decreases,
$s_0$ increases, and $\bar{s} \sim s_0$ as is easily seen by expanding
(\ref{soln1s}) in this limit.
Third, we remark that M.~Wolf obtains a similar result~\cite{newwolf},
except that he assumes {\it a priori} that $1\big/\bar{s} << 1$ in
order to simplify his analysis whereas our result is exact.

\subsection{The General Case}

Formally performing the summations and rearranging, (\ref{norm}), 
(\ref{snaught}), and (\ref{cut-off}) may be cast into the following 
useful forms (letting $q$ denote $\exp(-1\big/\bar{s})$):
\begin{equation}
\begin{split}
1 + \frac{f}{\pi_2} &= \frac{A}{1 - q} \, , \cr
1 &= \frac{A}{1-q} \, \left( 1 - q^{L+1} \right) \, , \cr
s_0 &= \frac{q}{1-q} - (L + 1) \, \frac{f}{\pi_2} \, . \cr
\end{split}
\end{equation}
These may be solved numerically for the set of parameters: $A$, $\bar{s}$
and $L$ once the inputs $s_0(\pi_{1} , \pi_{2})$, $\pi_2$, and $f$ are
specified.

Instead, to get an idea of the general behaviour, we make a well-motivated
approximation to the final equation listed.
In the regime described by the model, $f << \pi_2$, for reasonable values of 
$f$, and we expect that the factor of $L + 1$ is insufficient to render the
second term on the $rhs$ appreciable.
Put another way, the average value of the separation is rather insensitive
to the precise value of the maximum separation because relatively few 
separation events are maximal or near maximal.
Writing
\begin{equation}
s_0 = \frac{q}{1-q} \, ,
\label{s0approx}
\end{equation}
we immediately see that the solution for $\bar{s}$ is exactly the same as
above (\ref{soln1s}).
While this is an approximation, it is also a consequence of the relative
insensitivity to the cut-off as was strongly suggested by the success 
of the continuous analysis in \cite{random}.
With this result in hand, it is possible to solve for $A$ and $L$
without further approximations.
The normalisation constant
\begin{equation}
A = \left( 1 + \frac{f}{\pi_2} \right) \, \frac{1}{1 + s_0} 
\label{soln2A}
\end{equation}
is shifted slightly greater than its value derived in the case of no cut-off,
while
\begin{equation}
L = -1 + \frac{ \log\left( 1 + \frac{ \pi_2 }{f} \right) }%
{ \log\left[ 1 + \frac{1}{s_0} \right] }  \, ,
\label{soln2L}
\end{equation}
describes the growth in the cutoff in terms of the prime counts, statistical
average separation, and risk factor.

In Figure~\ref{graph3} we have replotted the actual thresholds 
obtained in our analysis of the likely maximal separations
as performed in \cite{return}.
Furthermore, we now include values of $L$ predicted
by (\ref{soln2L}) using Nicely's counts of primes and twins as inputs,
and choosing the risk factor to be equal to $1$.

{\begin{figure}[htb]
\begin{center}
\begin{turn}{-90}
\leavevmode
\epsfxsize=3.5in
\epsfysize=4.5in
\epsfbox{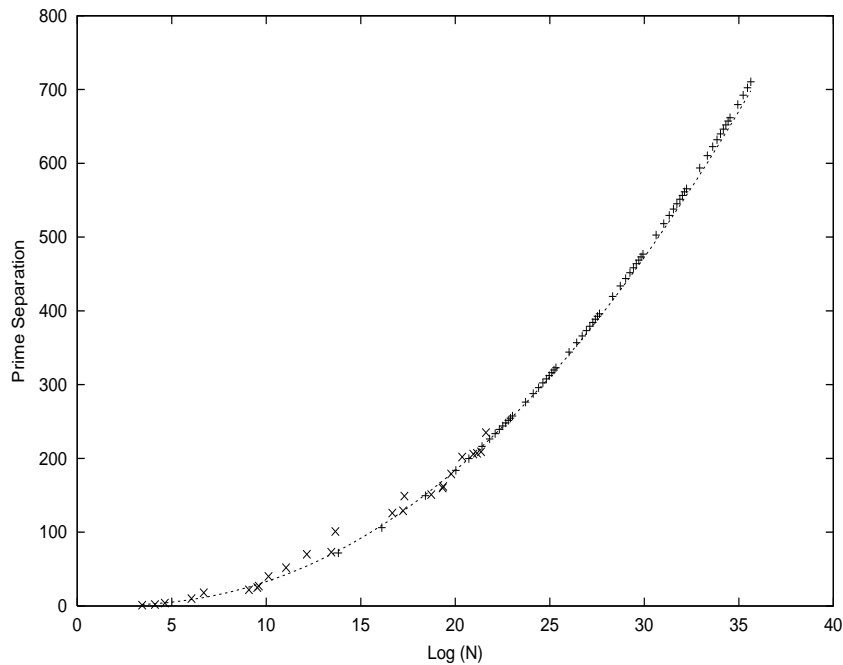}
\end{turn}
\end{center}
\caption{Maximal expected Prime Separation {\it vs.} $\log(N)$.
The $\times$'s in the lower left mark the actual onsets of
successive maxima in the gap spectrum.  The dotted curve is the
prediction for $f=1$ from the analysis in~\cite{return}.
The {\tt +}'s denote the computed values of $L$ from (\ref{soln2L}).}
\label{graph3}
\end{figure}}
 
In the region in which there is overlap between the observed thresholds,
our earlier computed estimates for the maximum separation at given $N$,
and the results obtained by direct computation using (\ref{soln2L}) 
the agreement is exceptional.
The spectacular agreement between the predicted cut-offs and
the maximal gap predictions of~\cite{return} persists above the range 
over which we have data.
While this was not unexpected, it provides further reassurance of the 
consistency of our model for the distribution of twins.

\section{Conclusion}

Very careful examination of the curve in Figure~\ref{graph2} reveals
that the data exhibit a slight tendency suggestive of negative curvature.
Naive theoretical considerations suggest attempting a three-parameter
fit of the form
\begin{equation}
s_{0}( \pi_1 ) = {\cal S}_{2} \log(\log( \pi_{1} ) ) + 
{\cal S}_{1} \log( \pi_{1} ) + {\cal S}_{0} \, ,
\label{sbar2}
\end{equation}
Fitting to the same data as before the empirical values obtained for
the constants ${\cal S}_{0}$, ${\cal S}_{1}$, and ${\cal S}_{2}$ are
\begin{equation*}
{\cal S}_{0} = -3.55 \pm 0.07 \ , \ 
{\cal S}_{1} = 0.745 \pm 0.001 \ , \mbox{ and }
{\cal S}_{2} = 1.10 \pm 0.03 \, .
\end{equation*}
From a practical perspective, it will require considerable effort to
extend the data into a regime in which the three parameter fit is 
clearly distinguished from the linear approximation.
This is just a manifestation of the extremely slow growth of the
function $\log_{2}(x) = \log( \log(x) )$.
We also note that the pathology of negative intercepts reappears with
increased strength.
Again, however, we can claim that this pathology is benign since it
implies a lower limit $\pi_{1}(N_{0}) \simeq \exp( 4.5 ) \sim 100$,
leading to an estimate that $N_0$ is on the order of $500$, which is
also far below the regime in which our statistical model is applicable.
Furthermore, the precise values obtained for ${\cal S}_{0}$ and
${\cal S}_{2}$ were rather sensitive to the range of data over which 
the fit was performed which leads us to believe that safer estimates 
for these coefficients and their errors are
\begin{equation*}
{\cal S}_{0} = -3.3 \pm 0.4 \ , \ 
{\cal S}_{1} = 0.75 \pm 0.01 \ , \mbox{ and }
{\cal S}_{2} = 1.0 \pm 0.2 \, .
\end{equation*}

In this paper, we have strengthened the case for, and extended the 
utility of, our characterisation of the the distribution of twin primes 
as random fixed-probability ``events'' among the primes.
We have done this by first performing the analysis of the model
without passing to the continuous (integral) limit and have 
demonstrated that there are no obstacles.
Furthermore we have made new and self-consistent predictions for the 
occurrence of ``gaps'' (maximum separations) and these are seen to 
conform well to the actual data and to our previous model analysis.
Perhaps the most important result going forward is our successful
reformulation of the model which has enabled it to accept as inputs
the raw counts of primes and twins below a certain $N$, along with
a risk-factor (of order 1) to which the model is fairly insensitive.

Lastly, we note that, particularly in this reformulated form without
accounting for the cutoff, the empirical model is {\it predictive}
since the factors which enter into the {\it Ansatz} (\ref{sprob}),
$A$ and $\bar{s}$, are determined from
$s_0$ alone whose behaviour is captured by (\ref{snaughtempirical}).
Given additional knowledge of $\pi_{2}(N)$ we can choose a risk-factor and
determine a more precise prediction for the spectrum including maximal
expected prime separation.

\section{Acknowledgements}

PFK and TP thank J.~Calvo and J.~Coykendall for helpful comments.
This work was supported in part by the National Science Foundation
(USA) under grant \#OSR-9452892 and an NSF EPSCoR Doctoral Dissertation 
Fellowship.


\begin{thebibliography}{50}
%
\bibitem{random} P.F.~Kelly and Terry Pilling, {\it Characterization of the
Distribution of Twin Primes}, math.NT/0103191.
\bibitem{return} P.F.~Kelly and Terry Pilling, {\it Implications of a
New Model of the Distribution of Twin Primes}, math.NT/0104205.
\bibitem{newwolf} Marek Wolf, {\it Some Remarks on the Distribution of Twin Primes}, math.NT/0105211; private communication.
\bibitem{jorge} J.~Calvo, private communication.
\bibitem{nicelydata} Thomas R.~Nicely, The tabulated values of $\pi_1(N)$ and $\pi_2(N)$ can be found at the website ``http://www.trnicely.net/index.html''

\end{thebibliography}
\end{document}